\documentstyle{amsppt}
\voffset-10mm
\magnification1200
\pagewidth{130mm}
\pageheight{204mm}
\hfuzz=2.5pt\rightskip=0pt plus1pt
\binoppenalty=10000\relpenalty=10000\relax
\TagsOnRight
\loadbold
\nologo
\addto\tenpoint{\normalbaselineskip=1.05\normalbaselineskip\normalbaselines}
\let\le\leqslant
\let\ge\geqslant

\topmatter
\title
On a combinatorial problem of Asmus Schmidt
\endtitle
\author
Wadim Zudilin\footnotemark"$^\ddag$"\ {\rm(Cologne \& Moscow)}
\endauthor
\date
\hbox to100mm{\vbox{\hsize=100mm%
\centerline{E-print \tt math.CA/0311195}
\smallskip
\centerline{August 2003}
}}
\enddate
\address
\hbox to70mm{\vbox{\hsize=70mm%
\leftline{Moscow Lomonosov State University}
\leftline{Department of Mechanics and Mathematics}
\leftline{Vorobiovy Gory, GSP-2, Moscow 119992 RUSSIA}
\leftline{{\it URL\/}: \tt http://wain.mi.ras.ru/index.html}
}}
\endaddress
\email
{\tt wadim\@ips.ras.ru}
\endemail
\abstract
For any integer $r\ge2$, define a sequence of
numbers $\{c_k^{(r)}\}_{k=0,1,\dots}$, independent
of the parameter~$n$, by
$$
\sum_{k=0}^n\binom nk^r\binom{n+k}k^r
=\sum_{k=0}^n\binom nk\binom{n+k}kc_k^{(r)},
\qquad n=0,1,2,\dotsc.
$$
We prove that all the numbers $c_k^{(r)}$ are integers.
\endabstract
\endtopmatter
\leftheadtext{W.~Zudilin}
\rightheadtext{On a problem of Asmus Schmidt}
\footnotetext"$^\ddag$"{The work is supported by an Alexander von Humboldt
research fellowship and partially supported by grant
no.~03-01-00359 of the Russian Foundation for Basic Research.}
\document

\head
1. Stating the problem
\endhead

The following curious problem was stated by A.\,L.~Schmidt
in~\cite{Sc1} in~1992.

\proclaim{Problem 1}
For any integer $r\ge2$, define a sequence of
numbers $\{c_k^{(r)}\}_{k=0,1,\dots}$, independent
of the parameter~$n$, by
$$
\sum_{k=0}^n\binom nk^r\binom{n+k}k^r
=\sum_{k=0}^n\binom nk\binom{n+k}kc_k^{(r)},
\qquad n=0,1,2,\dotsc.
\tag1
$$
Is it then true that all the numbers $c_k^{(r)}$ are integers\rom?
\endproclaim

An affirmative answer for $r=2$ was given in 1992
(but published a little bit later), independently,
by Schmidt himself~\cite{Sc2} and by V.~Strehl~\cite{St}.
They both proved the following explicit expression:
$$
c_n^{(2)}
=\sum_{j=0}^n\binom nj^3
=\sum_j\binom nj^2\binom{2j}n,
\qquad n=0,1,2,\dots,
\tag2
$$
which was observed experimentally by W.~Deuber, W.~Thumser
and B.~Voigt. In fact, Strehl used in~\cite{St} the corresponding
identity as a model for demonstrating various proof techniques
of binomial identities. He also proved an explicit
expression for the sequence $c_n^{(3)}$, thus answering affirmatively
to Problem~1 in the case $r=3$. But for this case Strehl had
only one proof based on Zeilberger's algorithm of creative telescoping.
Problem~1 was restated in~\cite{GKP}
(the last Research Problem on p.~256) with indication (on p.~549)
that H.\,S.~Wolf had shown the desired integrality of $c_n^{(r)}$
for any $r$ but only for any $n\le9$.

We recall that the first non-trivial
case $r=2$ is deeply related to
the famous Ap\'ery numbers $\sum_k\binom nk^2\binom{n+k}k^2$,
the denominators of rational approximations
to $\zeta(3)$. These numbers satisfy a 2nd-order polynomial
recursion discovered by R.~Ap\'ery in~1978, while an analogous
recursion (also 2nd-order and polynomial) for the numbers~\thetag{2}
was indicated by J.~Franel already in~1894.

The aim of this paper is to give an answer in the affirmative to Problem~1
(Theorem~1) by deriving explicit expressions for the numbers
$c_n^{(r)}$, and also to prove a stronger result (Theorem~2)
conjectured in~\cite{St, Section~4.2}.

\proclaim{Theorem 1}
The answer to Problem~\rom1 is affirmative. In particular,
we have the explicit expressions
$$
\align
c_n^{(4)}
&=\sum_j\binom{2j}j^3\binom nj
\sum_k\binom{k+j}{k-j}\binom j{n-k}\binom kj\binom{2j}{k-j},
\tag3
\\
c_n^{(5)}
&=\sum_j\binom{2j}j^4\binom nj^2
\sum_k\binom{k+j}{k-j}^2\binom{2j}{n-k}\binom{2j}{k-j},
\tag4
\endalign
$$
and in general for $s=1,2,\dots$
$$
\align
c_n^{(2s)}
&=\sum_j\binom{2j}j^{2s-1}\binom nj
\sum_{k_1}\binom j{n-k_1}\binom{k_1}j\binom{k_1+j}{k_1-j}
\sum_{k_2}\binom{2j}{k_1-k_2}\binom{k_2+j}{k_2-j}^2\dotsb
\\ &\qquad\times
\sum_{k_{s-1}}\binom{2j}{k_{s-2}-k_{s-1}}\binom{k_{s-1}+j}{k_{s-1}-j}^2
\binom{2j}{k_{s-1}-j},
\\
c_n^{(2s+1)}
&=\sum_j\binom{2j}j^{2s}\binom nj^2
\sum_{k_1}\binom{2j}{n-k_1}\binom{k_1+j}{k_1-j}^2
\sum_{k_2}\binom{2j}{k_1-k_2}\binom{k_2+j}{k_2-j}^2\dotsb
\\ &\qquad\times
\sum_{k_{s-1}}\binom{2j}{k_{s-2}-k_{s-1}}\binom{k_{s-1}+j}{k_{s-1}-j}^2
\binom{2j}{k_{s-1}-j},
\endalign
$$
where $n=0,1,2,\dots$\,.
\endproclaim

\head
2. Very-well-poised preliminaries
\endhead

The right-hand side of~\thetag{1} defines the so-called
{\it Legendre transform\/} of the sequence $\{c_k^{(r)}\}_{k=0,1,\dots}$.
In general, if
$$
a_n=\sum_{k=0}^n\binom nk\binom{n+k}kc_k
=\sum_{k=0}^n\binom{2k}k\binom{n+k}{n-k}c_k,
$$
then by the well-known relation for inverse Legendre pairs one has
$$
\binom{2n}nc_n=\sum_k(-1)^{n-k}d_{n,k}a_k,
$$
where
$$
d_{n,k}=\binom{2n}{n-k}-\binom{2n}{n-k-1}
=\frac{2k+1}{n+k+1}\binom{2n}{n-k}.
$$
Therefore, putting
$$
t_{n,j}^{(r)}=\sum_{k=j}^n(-1)^{n-k}d_{n,k}\binom{k+j}{k-j}^r,
\tag5
$$
we obtain
$$
\binom{2n}nc_n^{(r)}
=\sum_{j=0}^n\binom{2j}j^rt_{n,j}^{(r)}.
\tag6
$$
The case $r=1$ of Problem~1 is trivial (that is why it is not
included in the statement of the problem), while the cases
$r=2$ and $r=3$ are treated in~\cite{Sc2},~\cite{St} using the
fact that $t_{n,j}^{(2)}$ and $t_{n,j}^{(3)}$ have a {\it closed form}.
Namely, it is easy to show by Zeilberger's algorithm of
creative telescoping~\cite{PWZ} that the latter sequences, indexed by either
$n$ or~$j$, satisfy simple 1st-order polynomial recursions.
Unfortunately, this argument does not exist for $r\ge4$.

V.~Strehl observed in~\cite{St, Section~4.2} that the desired
integrality would be a consequence of the divisibility
of the product $\binom{2j}j^r\cdot t_{n,j}^{(r)}$ by $\binom{2n}n$
for all $j$, $0\le j\le n$. He conjectured a much stronger
property, which we are now able to prove.

\proclaim{Theorem 2}
The numbers $\binom{2n}n^{-1}\binom{2j}jt_{n,j}^{(r)}$
are integers.
\endproclaim

Our general strategy of proving Theorem~2 (and hence Theorem~1)
is as follows: rewrite \thetag{5} in a hypergeometric form and apply
suitable summation and transformation formulae (Propositions~1 and~2
below).

Changing $l$ to $n-k$ in~\thetag{5} we obtain
$$
t_{n,j}^{(r)}=\sum_{l\ge0}(-1)^l
\frac{2n-2l+1}{2n-l+1}\binom{2n}l\binom{n-l+j}{n-l-j}^r,
$$
where the series on the right terminates.
It is convenient to write all such terminating sums simply as $\sum_l$,
which is, in fact, a standard convention (see, e.g., \cite{PWZ}).
The ratio of the two consecutive terms in the latter sum is equal to
$$
\frac{-(2n+1)+l}{1+l}
\cdot\frac{-\frac12(2n-1)+l}{-\frac12(2n+1)+l}
\cdot\biggl(\frac{-(n-j)+l}{-(n+j)+l}\biggr)^r,
$$
hence
$$
t_{n,j}^{(r)}
=\binom{n+j}{n-j}^r
\cdot{}_{r+2}F_{r+1}\biggl(\matrix\format&\,\c\\
-(2n+1), & -\tfrac12(2n-1), & -(n-j), & \dots, & -(n-j) \\
& -\tfrac12(2n+1), & -(n+j), & \dots, & -(n+j)
\endmatrix\biggm|1\biggr)
$$
is a very-well-poised hypergeometric series. (We will omit
the argument $z=1$ in further discussions.)

The following two classical results---Dougall's summation
of a ${}_5F_4(1)$-series (proved in~1907) and
Whipple's transformation of a ${}_7F_6(1)$-series
(proved in~1926)---will be required to treat the cases $r=3,4,5$
of Theorems~1 and~2.

\proclaim{Proposition 1 \cite{Ba, Section~4.3}}
We have
$$
{}_5F_4\biggl(\matrix\format&\,\c\\
a, & 1+\frac12a, & c, & d, & -m \\
& \frac12a, & 1+a-c, & 1+a-d, & 1+a+m
\endmatrix\biggr)
=\frac{(1+a)_m\,(1+a-c-d)_m}{(1+a-c)_m\,(1+a-d)_m}
\tag7
$$
and
$$
\align
&
{}_7F_6\biggl(\matrix\format&\,\c\\
a, & 1+\frac12a, & b, & c, & d, & e, & -m \\
& \frac12a, & 1+a-b, & 1+a-c, & 1+a-d, & 1+a-e, & 1+a+m
\endmatrix\biggr)
\\ &\qquad
=\frac{(1+a)_m\,(1+a-d-e)_m}{(1+a-d)_m\,(1+a-e)_m}
\cdot{}_4F_3\biggl(\matrix
1+a-b-c, \, d, \, e, \, -m \\
1+a-b, \, 1+a-c, \, d+e-a-m
\endmatrix\biggr),
\tag8
\endalign
$$
where $m$ is a non-negative integer.
\endproclaim

Application of~\thetag{7} gives (without creative telescoping)
$$
t_{n,j}^{(3)}
=\binom{n+j}{n-j}^3
\cdot\frac{(-2n)_{n-j}(-2n+2(n-j))_{n-j}}{(-2n+(n-j))_{n-j}^2}
=\frac{(2n)!}{(3j-n)!\,(n-j)!^3},
$$
which is exactly the expression obtained in \cite{St, Section~4.2}.
Therefore, from~\thetag{6} we have the explicit expression
$$
c_n^{(3)}
=\binom{2n}n^{-1}\sum_j\binom{2j}j^3
\frac{(2n)!}{(3j-n)!\,(n-j)!^3}
=\sum_j\binom{2j}j^2\binom{2j}{n-j}\binom nj^2.
$$

For the case $r=5$, we are able to apply the transformation~\thetag{8}:
$$
\align
t_{n,j}^{(5)}
&=\binom{n+j}{n-j}^5
\cdot\frac{(-2n)_{n-j}(-2n+2(n-j))_{n-j}}{(-2n+(n-j))_{n-j}^2}
\\ &\qquad\times 
{}_4F_3\biggl(\matrix
-2j, \, -(n-j), \, -(n-j), \, -(n-j) \\
-(n+j), \, -(n+j), \, 3j-n+1
\endmatrix\biggr)
\\
&=\binom{n+j}{n-j}^2\frac{(2n)!}{(3j-n)!\,(n-j)!^3}
\sum_l\frac{(-2j)_l\,(-(n-j))_l^3}{l!\,(-(n+j))_l^2(3j-n+1)_l}
\\
&=\frac{(2n)!}{(2j)!\,(n-j)!^2}
\sum_l\binom{n-l+j}{n-l-j}^2\binom{2j}l
\binom{2j}{n-l-j}
\\
&=\frac{(2n)!}{(2j)!\,(n-j)!^2}
\sum_k\binom{k+j}{k-j}^2\binom{2j}{n-k}
\binom{2j}{k-j},
\endalign
$$
hence
$$
\binom{2n}n^{-1}\binom{2j}jt_{n,j}^{(5)}
=\binom nj^2
\sum_k\binom{k+j}{k-j}^2\binom{2j}{n-k}
\binom{2j}{k-j}
$$
are integers and from~\thetag{6} we derive formula~\thetag{4}.

To proceed in the case $r=4$, we apply the version of
formula~\thetag{8} with $b=(1+a)/2$ (so that the series
on the left reduces to a ${}_6F_5(1)$-very-well-poised series):
$$
\align
t_{n,j}^{(4)}
&=\binom{n+j}{n-j}^4
\cdot\frac{(-2n)_{n-j}(-2n+2(n-j))_{n-j}}{(-2n+(n-j))_{n-j}^2}
\\ &\qquad\times 
{}_4F_3\biggl(\matrix
-j, \, -(n-j), \, -(n-j), \, -(n-j) \\
-n, \, -(n+j), \, 3j-n+1
\endmatrix\biggr)
\allowdisplaybreak
&=\binom{n+j}{n-j}\frac{(2n)!}{(3j-n)!\,(n-j)!^3}
\sum_l\frac{(-j)_l\,(-(n-j))_l^3}{l!\,(-n)_l\,(-(n+j))_l(3j-n+1)_l}
\\
&=\frac{(2n)!\,j!}{n!\,(n-j)!\,(2j)!}
\sum_l\binom{n-l+j}{n-l-j}\binom jl\binom{n-l}j\binom{2j}{n-l-j}
\\
&=\frac{(2n)!\,j!}{n!\,(n-j)!\,(2j)!}
\sum_k\binom{k+j}{k-j}\binom j{n-k}\binom kj\binom{2j}{k-j},
\endalign
$$
from which, again, $\binom{2n}n^{-1}\binom{2j}jt_{n,j}^{(4)}\in\Bbb Z$
and we arrive at formula~\thetag{3}.

\head
3. Andrews's multiple transformation
\endhead

It seems that `classical' hypergeometric identities can
cover only the cases\footnote{%
This is not really true since Andrews's `non-classical' identity below
is a consequence of very classical Whipple's transformation and
the Pfaff--Saalsch\"utz formula.}
$r=2,3,\allowmathbreak4,5$ of Theorems~1 and~2.
In order to prove the theorems in full generality,
we will require a multiple generalization of Whipple's
transformation~\thetag{8}. The required generalization
is given by G.\,E.~Andrews in~\cite{An, Theorem~4}.
After making the passage $q\to1$ in Andrews's theorem,
we arrive at the following result.

\proclaim{Proposition 2}
For $s\ge1$ and $m$ a non-negative integer,
$$
\align
&
{}_{2s+3}F_{2s+2}\biggl(\matrix\format&\,\c\\
a, & 1+\frac12a, & b_1, & c_1, & b_2, & c_2, & \dots \\
& \frac12a, & 1+a-b_1, & 1+a-c_1, & 1+a-b_2, & 1+a-c_2, & \dots
\endmatrix
\\ &\qquad\qquad\qquad\qquad\qquad\qquad\qquad
\matrix\format&\,\c\\
\dots, & b_s, & c_s, & -m \\
\dots, & 1+a-b_s, & 1+a-c_s, & 1+a+m
\endmatrix\biggr)
\\ &\quad
=\frac{(1+a)_m(1+a-b_s-c_s)_m}{(1+a-b_s)_m(1+a-c_s)_m}
\sum_{l_1\ge0}\frac{(1+a-b_1-c_1)_{l_1}(b_2)_{l_1}(c_2)_{l_1}}
{l_1!\,(1+a-b_1)_{l_1}(1+a-c_1)_{l_1}}
\\ &\quad\qquad\times
\sum_{l_2\ge0}\frac{(1+a-b_2-c_2)_{l_2}(b_3)_{l_1+l_2}(c_3)_{l_1+l_2}}
{l_2!\,(1+a-b_2)_{l_1+l_2}(1+a-c_2)_{l_1+l_2}}
\dotsb
\\ &\quad\qquad\times
\sum_{l_{s-1}\ge0}
\frac{(1+a-b_{s-1}-c_{s-1})_{l_{s-1}}
(b_s)_{l_1+\dots+l_{s-1}}(c_s)_{l_1+\dots+l_{s-1}}}
{l_{s-1}!\,(1+a-b_{s-1})_{l_1+\dots+l_{s-1}}
(1+a-c_{s-1})_{l_1+\dots+l_{s-1}}}
\\ &\quad\qquad\qquad\times
\frac{(-m)_{l_1+\dots+l_{s-1}}}
{(b_s+c_s-a-m)_{l_1+\dots+l_{s-1}}}.
\endalign
$$
\endproclaim

\demo{Proof of Theorem~\rom2}
As in Section~2, we will distinguish the cases corresponding
to the parity of~$r$.

If $r=2s+1$, then setting $a=-(2n+1)$ and
$b_1=c_1=\dots=b_s=a_s=-m=\allowmathbreak-(n-j)$
in Proposition~2 we obtain
$$
\align
t_{n,j}^{(2s+1)}
&=\binom{n+j}{n-j}^{2s-2}\frac{(2n)!}{(3j-n)!\,(n-j)!^3}
\sum_{l_1}\binom{2j}{l_1}
\biggl(\frac{(-(n-j))_{l_1}}{(-(n+j))_{l_1}}\biggr)^2
\\ &\qquad\times
\sum_{l_2}\binom{2j}{l_2}
\biggl(\frac{(-(n-j))_{l_1+l_2}}{(-(n+j))_{l_1+l_2}}\biggr)^2
\dotsb
\\ &\qquad\times
\sum_{l_{s-1}}\binom{2j}{l_{s-1}}
\biggl(\frac{(-(n-j))_{l_1+\dots+l_{s-1}}}
{(-(n+j))_{l_1+\dots+l_{s-1}}}\biggr)^2
\\ &\qquad\qquad\times
\frac{(-1)^{l_1+\dots+l_{s-1}}(-(n-j))_{l_1+\dots+l_{s-1}}}
{(3j-n+1)_{l_1+\dots+l_{s-1}}}
\allowdisplaybreak
&=\frac{(2n)!}{(2j)!\,(n-j)!^2}
\sum_{l_1}\binom{2j}{l_1}\binom{n-l_1+j}{n-l_1-j}^2
\sum_{l_2}\binom{2j}{l_2}\binom{n-l_1-l_2+j}{n-l_1-l_2-j}^2
\dotsb
\\ &\qquad\times
\sum_{l_{s-1}}\binom{2j}{l_{s-1}}
\binom{n-l_1-\dotsb-l_{s-1}+j}{n-l_1-\dotsb-l_{s-1}-j}^2
\cdot\binom{2j}{n-l_1-\dotsb-l_{s-1}-j}.
\endalign
$$

If $r=2s$, we apply Proposition~2 with the choice
$a=-(2n+1)$, $b_1=(a+1)/2=-n$ and
$c_1=b_2=\dots=b_s=a_s=-m=-(n-j)$:
$$
\align
t_{n,j}^{(2s)}
&=\binom{n+j}{n-j}^{2s-3}\frac{(2n)!}{(3j-n)!\,(n-j)!^3}
\sum_{l_1}\binom j{l_1}\frac{(-(n-j))_{l_1}}{(-n)_{l_1}}
\,\frac{(-(n-j))_{l_1}}{(-(n+j))_{l_1}}
\\ &\qquad\times
\sum_{l_2}\binom{2j}{l_2}
\biggl(\frac{(-(n-j))_{l_1+l_2}}{(-(n+j))_{l_1+l_2}}\biggr)^2
\dotsb
\\ &\qquad\times
\sum_{l_{s-1}}\binom{2j}{l_{s-1}}
\biggl(\frac{(-(n-j))_{l_1+\dots+l_{s-1}}}
{(-(n+j))_{l_1+\dots+l_{s-1}}}\biggr)^2
\\ &\qquad\qquad\times
\frac{(-1)^{l_1+\dots+l_{s-1}}(-(n-j))_{l_1+\dots+l_{s-1}}}
{(3j-n+1)_{l_1+\dots+l_{s-1}}}
\allowdisplaybreak
&=\frac{(2n)!\,j!}{n!\,(n-j)!\,(2j)!}
\sum_{l_1}\binom j{l_1}\binom{n-l_1}j\binom{n-l_1+j}{n-l_1-j}
\\ &\qquad\times
\sum_{l_2}\binom{2j}{l_2}\binom{n-l_1-l_2+j}{n-l_1-l_2-j}^2
\dotsb
\\ &\qquad\times
\sum_{l_{s-1}}\binom{2j}{l_{s-1}}
\binom{n-l_1-\dotsb-l_{s-1}+j}{n-l_1-\dotsb-l_{s-1}-j}^2
\cdot\binom{2j}{n-l_1-\dotsb-l_{s-1}-j}.
\endalign
$$

In both cases, the desired integrality
$$
\binom{2n}n^{-1}\binom{2j}jt_{n,j}^{(r)}\in\Bbb Z,
\qquad j=0,1,\dots,n,
$$
clearly holds, and Theorem~2 follows.
\enddemo

Theorem~1 is an immediate consequence of Theorem~2.

\medskip
We would like to conclude the paper by the following $q$-question.

\proclaim{Problem 2}
Find and solve an appropriate $q$-analogue of Problem~\rom1.
\endproclaim

\subsubhead
Acknowledgements
\endsubsubhead
I was greatly encouraged by C.~Krattenthaler to prove
binomial identities by myself. I thank him for our
fruitful discussions and for pointing out to me Andrews's formula.
I thank J.~Sondow for several suggestions that allowed
me to improve the text of the paper.
This work was done during a long-term visit
at the Mathematical Institute of Cologne University.
I thank the staff of the institute and personally
P.~Bundschuh for the brilliant working atmosphere I had there.


\enddocument